\documentclass{ws-procs975x65}

\def\beq{\begin{equation}}
\def\eeq{\end{equation}}

\begin{document}

\title{A non-geometric representation of the Dirac equation\\ in curved spacetime}
\author{Dmitri Vassiliev}

\address{Department of Mathematics,
University College London,\\ Gower Street, London WC1E 6BT, UK\\
E-mail: D.Vassiliev@ucl.ac.uk\\
www.homepages.ucl.ac.uk/$\sim$ucahdva/\\
Supported by EPSRC grant EP/M000079/1}

\begin{abstract}
We write the Dirac equation in curved 4-dimensional Lorentzian spacetime
using concepts from the analysis of partial differential equations as opposed
to geometric concepts.
\end{abstract}


\bodymatter


\section{Playing Field}

Let $M$ be a connected 4-manifold without boundary.
We will work with 2-columns $v:M\to\mathbb{C}^2$ of complex-valued
half-densities (a half-density is a quantity which transforms as the square root
of a density under changes of local coordinates).
The inner product on such 2-columns is defined as
$\langle v,w\rangle:=\int_M w^*v\,dx$,
where $x=(x^1,x^2,x^3,x^4)$ are local coordinates,
$dx=dx^1dx^2dx^3dx^4$
and the star stands for Hermitian conjugation.

Let $L$ be a formally self-adjoint first order linear differential
operator acting on 2-columns of complex-valued half-densities.
Our initial objective will be to examine the geometric content of the operator $L$.
In order to pursue this objective we first need to provide an invariant analytic
description of the operator.

In local coordinates our operator reads
\begin{equation}
\label{operator L in local coordinates}
L=F^\alpha(x)\frac\partial{\partial x^\alpha}+G(x),
\end{equation}
where $F^\alpha(x)$, $\alpha=1,2,3,4$, and $G(x)$ are some $2\times 2$ matrix-functions.
The principal and subprincipal symbols of the operator $L$ are defined as
\begin{equation}
\label{definition of the principal symbol}
L_\mathrm{prin}(x,p):=iF^\alpha(x)\,p_\alpha\,,
\end{equation}
\begin{equation}
\label{definition of the subprincipal symbol}
L_\mathrm{sub}(x):=G(x)
+\frac i2(L_\mathrm{prin})_{x^\alpha p_\alpha}(x)\,,
\end{equation}
where $p=(p_1,p_2,p_3,p_4)$ is the dual variable (momentum); see Ref.~\refcite{mybook}.
The principal and subprincipal symbols 
are invariantly defined $2\times 2$ Hermitian matrix-functions
on $T^*M$ and $M$ respectively
which uniquely determine the operator $L$.

Further on in this paper we assume that the principal symbol
of our operator satisfies the following
non-degeneracy condition:
\begin{equation}
\label{definition of non-degeneracy}
L_\mathrm{prin}(x,p)\ne0,\qquad\forall(x,p)\in T^*M\setminus\{0\}.
\end{equation}
Condition \eqref{definition of non-degeneracy}
means that the elements of the $2\times2$ matrix-function
$L_\mathrm{prin}(x,p)$ do not vanish simultaneously
for any $x\in M$ and any nonzero momentum $p$.

\section{Lorentzian Metric and Orthonormal Frame}
\label{Lorentzian Metric and Orthonormal Frame}

Observe that the determinant of the principal symbol is a quadratic form
in the dual variable (momentum) $p$\,:
\begin{equation}
\label{definition of metric}
\det L_\mathrm{prin}(x,p)=-g^{\alpha\beta}(x)\,p_\alpha p_\beta\,.
\end{equation}
We interpret the real coefficients $g^{\alpha\beta}(x)=g^{\beta\alpha}(x)$,
$\alpha,\beta=1,2,3,4$, appearing in formula \eqref{definition of metric}
as components of a (contravariant) metric tensor.

The following result was established in Ref.~\refcite{nongeometric}.

\begin{lemma}
\label{Lemma about Lorentzian metric}
Our metric is Lorentzian, i.e.~it has
three positive eigenvalues
and one negative eigenvalue.
\end{lemma}

Furthermore, the principal symbol of our operator defines an orthonormal
frame $e_j{}^\alpha(x)$.
Here the Latin index $j=1,2,3,4$ enumerates the vector fields,
the Greek index $\alpha=1,2,3,4$ enumerates the components
of a given vector $e_j$
and orthonormality is understood in the Lorentzian sense:
\begin{equation}
\label{orthonormality of the frame}
g_{\alpha\beta}\,e_j{}^\alpha e_k{}^\beta=
\begin{cases}
0\quad&\text{if}\quad j\ne k,
\\
1\quad&\text{if}\quad j=k\ne4,
\\
-1\quad&\text{if}\quad j=k=4.
\end{cases}
\end{equation}

The orthonormal frame is recovered from the principal symbol as follows.
Decomposing the principal symbol with respect to the standard basis
\begin{equation*}
\label{standard basis}
s^1=
\begin{pmatrix}
0&\,1\,\\
\,1\,&0
\end{pmatrix},
\quad
s^2=
\begin{pmatrix}
0&-i\\
i&0
\end{pmatrix},
\quad
s^3=
\begin{pmatrix}
1&0\\
0&-1
\end{pmatrix},
\quad
s^4=
\begin{pmatrix}
\,1\,&0\\
0&\,1\,
\end{pmatrix}
\end{equation*}
in the real vector space of $2\times2$ Hermitian matrices, we get
$L_\mathrm{prin}(x,p)=s^j c_j(x,p)$.
Each coefficient $c_j(x,p)$ is linear in momentum $p$, so $c_j(x,p)=e_j{}^\alpha(x)\,p_\alpha\,$.

The existence of an orthonormal frame implies that our manifold $M$ is parallelizable.
We see that our analytic non-degeneracy condition
\eqref{definition of non-degeneracy}
has far reaching geometric consequences.

\section{Gauge Transformations and Covariant Subprincipal Symbol}

Let us consider the action (variational functional)
$\,\int_M v^*(Lv)\,dx\,$
associated with our operator.
Take an arbitrary smooth matrix-function
\begin{equation}
\label{SL2C matrix-function R}
R:M\to\mathrm{SL}(2,\mathbb{C})
\end{equation}
and consider the following transformation of our 2-column of unknowns:
\begin{equation}
\label{SL2C transformation of the unknowns}
v\mapsto Rv.
\end{equation}
We interpret
\eqref{SL2C transformation of the unknowns}
as a gauge transformation
because we are looking here at a change of basis in our
vector space of unknowns $v:M\to\mathbb{C}^2$.

The transformation
\eqref{SL2C transformation of the unknowns}
of the 2-column $v$ induces the following transformation of the action:
$\,\int_M v^*(Lv)\,dx
\,\mapsto
\int_M v^*(R^*LRv)\,dx\,$.
This means that our $2\times2$ differential operator $L$ experiences the transformation
\begin{equation}
\label{SL2C transformation of the operator}
L\mapsto R^*LR\,.
\end{equation}
This section is dedicated to the analysis of the transformation
\eqref{SL2C transformation of the operator}.

\begin{remark}
We chose to restrict our analysis to matrix-functions $R(x)$ of determinant one,
see formula \eqref{SL2C matrix-function R},
because we want to preserve our Lorentzian metric defined in accordance
with formula \eqref{definition of metric}.
\end{remark}

\begin{remark}
In non-relativistic theory one normally looks at the transformation
\begin{equation}
\label{SU2 transformation of the operator}
L\mapsto R^{-1}LR
\end{equation}
rather than at
\eqref{SL2C transformation of the operator}.
The reason we chose to go along with \eqref{SL2C transformation of the operator}
is that we are thinking in terms of actions and corresponding
Euler--Lagrange equations rather than operators as such.
We believe that this point of view makes more sense in the relativistic setting.
If one were consistent in promoting such a point of view, then one would have had to
deal with actions throughout the paper rather than with operators.
We did not adopt this `consistent' approach because
this would have made the paper difficult to read.
Therefore, throughout the paper
we use the concept of an operator, having in mind that we are
really interested in the action and corresponding Euler--Lagrange equation.
\end{remark}

\begin{remark}
The transformations
\eqref{SL2C transformation of the operator}
and
\eqref{SU2 transformation of the operator}
coincide if the matrix-function $R(x)$ is special unitary.
Applying special unitary transformations is natural in the non-relativistic
3-dimensional setting when dealing with an elliptic system, see Ref.~\refcite{arxiv},
but in the relativistic
4-dimensional setting when dealing with a hyperbolic system
special unitary transformations are too restrictive.
\end{remark}


The transformation \eqref{SL2C transformation of the operator}
of the differential operator $L$ induces the following transformations
of its principal \eqref{definition of the principal symbol}
and subprincipal \eqref{definition of the subprincipal symbol}
symbols:
\begin{equation}
\label{SL2C transformation of the principal symbol}
L_\mathrm{prin}\mapsto R^*L_\mathrm{prin}R\,,
\end{equation}
\begin{equation}
\label{SL2C transformation of the subprincipal symbol}
L_\mathrm{sub}\mapsto
R^*L_\mathrm{sub}R
+\frac i2
\left(
R^*_{x^\alpha}(L_\mathrm{prin})_{p_\alpha}R
-
R^*(L_\mathrm{prin})_{p_\alpha}R_{x^\alpha}
\right).
\end{equation}
Comparing formulae
\eqref{SL2C transformation of the principal symbol}
and
\eqref{SL2C transformation of the subprincipal symbol}
we see that, unlike the principal symbol, the subprincipal
symbol does not transform in a covariant fashion due to
the appearance of terms with the gradient of the
matrix-function $R(x)$.

It turns out that one can overcome the non-covariance
in \eqref{SL2C transformation of the subprincipal symbol} by introducing
the \emph{covariant subprincipal symbol} $\,L_\mathrm{csub}(x)\,$
in accordance with formula
\begin{equation}
\label{definition of covariant subprincipal symbol}
L_\mathrm{csub}:=
L_\mathrm{sub}
+\frac i{16}\,
g_{\alpha\beta}
\{
L_\mathrm{prin}
,
\operatorname{adj}L_\mathrm{prin}
,
L_\mathrm{prin}
\}_{p_\alpha p_\beta},
\end{equation}
where
$
\{F,G,H\}:=F_{x^\alpha}GH_{p_\alpha}-F_{p_\alpha}GH_{x^\alpha}
$
is the generalised Poisson bracket on matrix-functions
and $\,\operatorname{adj}\,$ is the operator of matrix adjugation
\begin{equation}
\label{definition of adjugation}
F=\begin{pmatrix}a&b\\ c&d\end{pmatrix}
\mapsto
\begin{pmatrix}d&-b\\-c&a\end{pmatrix}
=:\operatorname{adj}F
\end{equation}
from elementary linear algebra.

The following result was established in Ref.~\refcite{nongeometric}.

\begin{lemma}
\label{Lemma about covariant subprincipal symbol}
The transformation
\eqref{SL2C transformation of the operator}
of the differential operator induces the transformation
$\,L_\mathrm{csub}\mapsto R^*L_\mathrm{csub}R\,$
of its covariant subprincipal symbol.
\end{lemma}

Comparing formulae
\eqref{definition of the subprincipal symbol}
and
\eqref{definition of covariant subprincipal symbol}
we see that the standard subprincipal symbol
and covariant subprincipal symbol have the same structure, only
the covariant subprincipal symbol has a second correction term
designed to `take care' of special linear transformations
in the vector space of unknowns $v:M\to\mathbb{C}^2$.
The standard subprincipal symbol \eqref{definition of the subprincipal symbol}
is invariant under changes of local coordinates
(its elements behave as scalars),
whereas the covariant subprincipal
symbol~\eqref{definition of covariant subprincipal symbol}
retains this feature but gains an extra $\mathrm{SL}(2,\mathbb{C})$
covariance property. In other words, the covariant subprincipal symbol
\eqref{definition of covariant subprincipal symbol}
behaves `nicely' under a wider group of transformations.

\section{Electromagnetic Covector Potential}
\label{Electromagnetic Covector Potential}

The covariant subprincipal
symbol can be uniquely represented in the form
\begin{equation}
\label{decomposition of  covariant subprincipal symbol}
L_\mathrm{csub}(x)=L_\mathrm{prin}(x,A(x)),
\end{equation}
where $A=(A_1,A_2,A_3,A_4)$ is some real-valued covector field.
We interpret this covector field as the electromagnetic covector potential.

Lemma~\ref{Lemma about covariant subprincipal symbol}
and formulae
\eqref{SL2C transformation of the principal symbol}
and
\eqref{decomposition of  covariant subprincipal symbol}
tell us that the electromagnetic covector potential
is invariant under gauge transformations
\eqref{SL2C transformation of the operator}.

\section{Adjugate Operator}

\begin{definition}
The adjugate of
a formally self-adjoint non-degenerate first order $2\times2$ linear  differential operator $L$
is
the formally self-adjoint non-degenerate first order $2\times2$ linear  differential operator
$\operatorname{Adj}L$
whose principal and covariant subprincipal symbols are matrix adjugates
of those of the operator $L$.
\end{definition}

We denote matrix adjugation by $\,\operatorname{adj}\,$,
see formula \eqref{definition of adjugation},
and operator adjugation by $\,\operatorname{Adj}\,$.
Of course,
the coefficients of the adjugate operator
can be written down explicitly in local coordinates via
the coefficients of the original operator \eqref{operator L in local coordinates},
see Ref.~\refcite{nongeometric} for details.

Applying the analysis from Sections
\ref{Lorentzian Metric and Orthonormal Frame}--\ref{Electromagnetic Covector Potential}
to the differential operator $\operatorname{Adj}L$ it is easy to see
that the metric and electromagnetic covector potential
encoded within the operator $\operatorname{Adj}L$ are the same as
in the original operator $L$.
Thus, the metric and electromagnetic covector potential
are invariant under operator adjugation.

It also easy to see that $\operatorname{Adj}\operatorname{Adj}L=L$,
so operator adjugation is an involution.

\section{Main Result}

We define the Dirac operator as the differential operator
\begin{equation}
\label{analytic definition of the Dirac operator}
D:=
\begin{pmatrix}
L&mI\\
mI&\operatorname{Adj}L
\end{pmatrix}
\end{equation}
acting on 4-columns
$\,\psi=\begin{pmatrix}
\,v_1\,&v_2\,&w_1\,&w_2\,
\end{pmatrix}^T\,$
of complex-valued half-densities.
Here $m$ is the electron mass and
$I$ is the $2\times2$ identity matrix.

The `traditional' Dirac operator $D_\mathrm{trad}$ is written down in Appendix A
of Ref.~\refcite{nongeometric} and acts on bispinor fields
$\,\psi_\mathrm{trad}=\begin{pmatrix}
\,\xi^1\,&\xi^2\,&\eta_{\dot 1}\,&\eta_{\dot 2}\,\end{pmatrix}^T\,$.
Here we assume, without loss of generality, that the orthonormal frame used in the
construction of the operator $D_\mathrm{trad}$ is the one from
Section~\ref{Lorentzian Metric and Orthonormal Frame}.

Our main result is the following theorem established in Ref.~\refcite{nongeometric}.

\begin{theorem}
\label{main theorem}
The two operators,
our analytically defined Dirac operator \eqref{analytic definition of the Dirac operator}
and geometrically defined Dirac operator $D_\mathrm{trad}\,$,
are related by the formula
\begin{equation}
\label{main theorem formula}
D
=
|\det g_{\kappa\lambda}|^{1/4}
\,
D_\mathrm{trad}
\,
|\det g_{\mu\nu}|^{-1/4}\,.
\end{equation}
\end{theorem}

Consider now the two Dirac equations
\begin{equation}
\label{our Dirac equation}
D\psi=0,
\end{equation}
\begin{equation}
\label{traditional Dirac equation}
D_\mathrm{trad}\psi_\mathrm{trad}=0.
\end{equation}
Formula \eqref{main theorem formula}
implies that the solutions of equations
\eqref{our Dirac equation}
and
\eqref{traditional Dirac equation}
differ only by a prescribed scaling factor:
$\,\psi=|\det g_{\mu\nu}|^{1/4}\,\psi_\mathrm{trad}\,$.
This means that for all practical purposes equations
\eqref{our Dirac equation}
and
\eqref{traditional Dirac equation}
are equivalent.

\section{Spin Structure}

Let us consider all possible formally self-adjoint non-degenerate first order $2\times2$
linear differential operators $L$
corresponding, in the sense of formula \eqref{definition of metric},
to the prescribed Lorentzian metric.
In this section our aim is to classify all such operators~$L$.

Let us fix a reference operator ${\mathbf L}$ and let
${\mathbf e}_j$ be the corresponding orthonormal frame
(see Section~\ref{Lorentzian Metric and Orthonormal Frame}).
Let $L$ be another operator and let
$e_j$ be the corresponding orthonormal frame.
We define the following two real-valued scalar fields
\[
\mathbf{c}(L):=
-\,\frac1{4!}\,(
{\mathbf e}_1
\wedge{\mathbf e}_2
\wedge{\mathbf e}_3
\wedge{\mathbf e}_4
)_{\kappa\lambda\mu\nu}
\,(e_1\wedge e_2\wedge e_3\wedge e_4)^{\kappa\lambda\mu\nu}\,,
\quad
\mathbf{t}(L):=
-\,{\mathbf e}_{4\alpha}\,e_4{}^\alpha\,.
\]
Observe that these scalar fields do not vanish; in fact,
$\mathbf{c}(L)$ can take only two values, $+1$ or $-1$.
This observation gives us a primary classification of operators $L$ into four classes
determined by the signs of $\mathbf{c}(L)$ and $\mathbf{t}(L)$.
The four classes correspond to the four connected components of
the Lorentz group.

Note that
\begin{eqnarray*}
&\mathbf{c}(-L)=\mathbf{c}(L),
\qquad
&\mathbf{t}(-L)=-\mathbf{t}(L),
\\
&\ \ \ \ \ \,\mathbf{c}(\operatorname{Adj}L)=-\mathbf{c}(L),
\qquad
&\mathbf{t}(\operatorname{Adj}L)=\mathbf{t}(L),
\end{eqnarray*}
which means that by applying the transformations $L\mapsto-L$ and
$L\mapsto\operatorname{Adj}L$
to a given operator $L$ one can reach all four classes of our
primary classification.

Further on we work with operators $L$ such that $\mathbf{c}(L)>0$ and $\mathbf{t}(L)>0$.

We say that the operators $L$ and $\tilde L$ are equivalent if there
exists a smooth matrix-function \eqref{SL2C matrix-function R}
such that $\tilde L_\mathrm{prin}=R^*L_\mathrm{prin}R$.
The equivalence classes of operators obtained this way
are called \emph{spin structures}.

The above 4-dimensional Lorentzian definition of spin structure is an extension of the
3-dimensional Riemannian definition from Ref.~\refcite{arxiv}.
The difference is that we have now dropped the condition
$\operatorname{tr}L_\mathrm{prin}(x,p)=0$,
replaced the ellipticity condition by the
weaker non-degeneracy condition \eqref{definition of non-degeneracy}
and extended our group of transformations
from special unitary to special linear.

One would hope that for a connected  Lorentzian 4-manifold
admitting  a global orthonormal frame
(see  \eqref{orthonormality of the frame} for definition of orthonormality)
our analytic definition of spin structure would be equivalent
to the traditional geometric one. Unfortunately, we do not currently
have a rigorous proof of equivalence in the 4-dimensional Lorentzian setting.

\section*{Acknowledgments}

The author is grateful to Zhirayr Avetisyan and Nikolai Saveliev for helpful advice.

\end{document}